\documentclass[preprint,12pt,authoryear]{elsarticle}

\usepackage{comment}
\usepackage{subfigure}
\usepackage{multirow}

\usepackage{algorithm}
\usepackage{algpseudocode}

\usepackage{hyperref}

\usepackage{natbib}

\usepackage{amssymb}
\usepackage{amsmath}
\allowdisplaybreaks

\usepackage{amssymb}
\usepackage{amsmath}

\journal{Computers \& Operations Research}

\begin{document}

\begin{frontmatter}



\title{Progressive hedging for multi-stage stochastic lot sizing problems with setup carry-over under uncertain demand}


\author[a]{Manuel Schlenkrich}
\author[b]{Jean-Fran\c{c}ois Cordeau}
\author[a]{Sophie N. Parragh}

\affiliation[a]{organization={Insitute of Production and Logistics Management, Johannes Kepler University Linz},
            addressline={Altenbergerstra\ss e 69}, 
            city={Linz},
            postcode={4040}, 
            country={Austria}}
            
\affiliation[b]{organization={Department of Logistics and Operations Management, HEC Montréal},
    addressline={\\3000 chemin de la Côte-Sainte-Catherine}, 
    city={Montréal},
    postcode={H3T 2A7}, 
    country={Canada}}

\begin{abstract}
We investigate multi-stage demand uncertainty for the multi-item multi-echelon capacitated lot sizing problem with setup carry-over. Considering a multi-stage decision framework helps to quantify the benefits of being able to adapt decisions to newly available information. The drawback is that multi-stage stochastic optimization approaches lead to very challenging formulations. This is because they usually rely on scenario tree representations of the uncertainty, which grow exponentially in the number of decision stages. Thus, even for a moderate number of decision stages it becomes difficult to solve the problem by means of a compact optimization model. To address this issue, we propose a progressive hedging algorithm and we investigate and tune the crucial penalty parameter that influences the conflicting goals of fast convergence and solution quality. While low penalty parameters usually lead to high quality solutions, this comes at the cost of slow convergence. To tackle this problem, we adapt metaheuristic adjustment strategies to guide the algorithm towards a consensus more efficiently. Furthermore, we consider several options to compute the consensus solution. While averaging the subproblem decisions is a common choice, we also apply a majority voting procedure. We test different algorithm configurations and compare the results of progressive hedging to the solutions obtained by solving a compact optimization model on well-known benchmark instances. For several problem instances the progressive hedging algorithm converges to solutions within 1\% of the cost of the compact model's solution, while requiring shorter runtimes.
\end{abstract}
\begin{highlights}
    \item Define a multi-stage stochastic lot sizing problem with setup carry-over
    \item Analyze the penalty parameter's role in the progressive hedging method's performance
    \item Apply metaheuristic strategies to enhance the algorithm's convergence behavior
    \item Find varying algorithmic behavior for different consensus calculation procedures
    \item Runtime of tuned progressive hedging outperforms compact model within 1\% cost range
\end{highlights}

\end{frontmatter}


\section{Introduction}
Lot sizing is a crucial component of production planning systems that aim to satisfy external demand in a cost-efficient manner. In addition, it has long been recognized that considering partial scheduling decisions in lot sizing leads to more realistic production plans and is especially important for problems involving long and costly production setups. Along these lines, lot sizing with setup carry-over (see, e.g., \cite{sox_capacitated_1999} and \cite{ghirardi_matheuristics_2019}) allows to preserve setups on production resources between consecutive periods. \cite{WEI201799} state that setup carry-overs represent a practically relevant aspect of production systems that operate continuously around the clock. These 24/7 production environments consider an uninterrupted planning horizon and make setup carry-overs necessary to provide the decision maker with the required level of flexibility \citep{belo-filho_models_2014}.

Lot sizing for capacitated multi-item multi-echelon production systems with setup carry-overs is a challenging task even in a deterministic setting. On top of that, customer demand is rarely known in advance but is often uncertain. This significantly impacts the respective lot sizing decisions. As \cite{sereshti_managing_2024} point out, a majority of the available literature focuses on problem settings that consider only end items to face external customer demand, while components are bound to internal requirements along the bill of materials (BOM). However, industries providing aftermarket services depend on spare parts and therefore the problem variant with external demand at multiple BOM levels is practically relevant \citep{wagner_lindemann_spare_parts}. For instance, in the aerospace industry there is not only a need for end items, but also components are required in the repair and maintenance processes. Planning processes must hence consider this external demand for spare parts alongside the demand for end items \citep{dodin_bombardier_2023}. 

Considering uncertain external demand for items across all BOM levels in a multi-stage stochastic setting leads to an additional level of decision, namely whether components should be instantly used to satisfy external customer demand or be further processed to satisfy end item demand in later periods.  

Stochastic demand has already been considered for the capacitated multi-item multi-echelon lot sizing problem in two-stage frameworks and also multi-stage decision processes have been recently addressed, such as in \cite{tunc_extended_2018} and \cite{thevenin_material_2021}. Uncertain customer demand was also incorporated into the integrated production planning and scheduling problem, a special variant considering detailed scheduling constraints, e.g.,~in \cite{hu_two-stage_2016} and \cite{alem_computational_2018}. In a stochastic context having detailed scheduling decisions in the second stage can lead to significant tractability issues. This problem becomes even more severe when considering a multi-stage stochastic setting. 
To tackle this issue, we consider the possibility of setup carry-overs, which allow to meet the necessary requirements of 24/7 production systems.

To the best of our knowledge, \cite{SCHLENKRICH2024109379} are among the first to study demand uncertainty for production environments with multiple production levels, limited resource capacities and the possibility to carry over production setups. While this work analyzes a two-stage approximation, we investigate the problem in a multi-stage setting and propose a progressive hedging algorithm to tackle this challenging problem. This helps to understand the benefits of being able to adapt decisions as soon as new information becomes available, especially in the presence of external demand for components. 

The remainder of the article is structured as follows. In Section \ref{sec:literature_review} we review the lot sizing literature with a focus on the aspects of setup carry-over and uncertainty. Section \ref{sec:problem_formulation_and_model} defines the studied lot sizing problem and a multi-stage stochastic optimization model is formulated. A progressive hedging solution approach is proposed in Section \ref{sec:progressive_hedging} and evaluated in an extensive numerical study in Section \ref{sec:computational_study}. Finally, in Section \ref{sec:conclusion} we conclude our work.

\section{Literature Review}
\label{sec:literature_review}

Several lot sizing surveys provide a general overview of the existing literature. \citet{karimi_capacitated_2003} review models and algorithms for the capacitated problem variant, while \citet{quadt_capacitated_2008} focus on lot sizing extensions. The survey of \citet{jans_modeling_2008} reviews industrial lot sizing problems and \citet{brahimi_single-item_2017} provide a recent overview for the single-item case. The remainder of our literature review is structured in two sections focusing on the problem characteristics considered in this work. We first discuss lot sizing with setup carry-over and then review uncertainty in lot sizing problems.

\subsection{Lot sizing with setup carry-over}
Lot sizing with \textit{setup carry-overs} also called \textit{linked production lots}, denotes a problem variant which allows to carry over the production setup of an item to the consecutive period. Hence, production of the respective items can continue in the following period, without the need to perform new setups. In the seminal work by \citet{sox_capacitated_1999}, the authors propose the concept of setup carry-over for production systems that allow production of multiple items in one period. The carry-over decisions are modeled by binary variables and we follow this approach in the present work.

After this seminal paper by \citet{sox_capacitated_1999}, the problem variant including setup carry-over was investigated by numerous authors. The value of considering setup carry-overs was demonstrated by \citet{porkka_multiperiod_2003}, who compare different lot sizing models. The authors show that models ignoring setup carry-overs are outperformed by models that allow carry-overs. Additional production planning aspects for the lot sizing variant with setup carry-overs were studied in numerous extensions. \citet{sung_mixed-integer_2008} study the case of very long setup procedures, in which setup times might exceed the length of a single planning period. The multi-plant lot sizing problem with setup carry-over is investigated by \citet{nascimento_hybrid_2008} and \citet{sahling_solving_2009} address a problem variant allowing for overtime production. Several authors have also proposed novel exact and heuristic solution methods for the lot sizing problem with setup carry-over, i.e.,~\citet{karimi_capacitated_2003} present a Tabu Search algorithm, \citet{briskorn_note_2006} use dynamic programming to enhance their heuristic approach and \citet{suerie_capacitated_2003} develop a decomposition method based on valid inequalities. 

In the literature, several works are motivated by practical applications of the lot sizing problem with setup carry-overs. For instance, \citet{quadt_capacitated_2009} study the case of semiconductor production, which usually requires the use of several parallel machines. In the work of \citet{ghirardi_matheuristics_2019}, additionally the respective setup and processing times vary among the different machines. The special variant of so-called setup cross-overs, presented by \citet{belo-filho_models_2014}, provides the opportunity to continue the setup procedure in the following period. \citet{gansterer_capacitated_2021} consider the possibility of transshipment, as well as setup carry-over, in order to allow collaboration between multiple agents.

In the literature, only a few papers consider the lot sizing problem with setup carry-over including uncertainty.
\citet{behnamian_markovian_2017} consider a production setting in which products are inspected regarding their quality and might require additional production steps or are sorted out as scrap. The work therefore investigates uncertainty regarding the production levels, while assuming deterministic demand. In a recent work, \citet{SCHLENKRICH2024109379} propose robust as well as two-stage stochastic lot sizing models for multi-level production systems under demand uncertainty and compare the different approaches regarding their suitability for various problem characteristics. Lastly, \cite{SIMONIS2025} present a simulation-optimization approach for tablet production in the pharmaceutical industry. The authors consider uncertain customer demand and present a two-stage stochastic lot sizing model with setup carry-over and compare it to a generalized uncertainty framework. Even though the problem variant has already been considered in a two-stage stochastic setting, the adaptive nature of the underlying multi-stage decision process has not been investigated.

\subsection{Lot sizing under uncertainty}
Multi-stage stochastic optimization problems usually lead to extremely large mathematical models, which are challenging to solve. It is thus often beneficial to tackle these sequential decision problems by applying (approximate) dynamic programming and reinforcement learning, but also decomposition methods, such as progressive hedging, as well as heuristics. 

Dynamic programming and reinforcement learning represent powerful techniques to derive policies for stochastic sequential decision processes and have already been applied in the lot sizing context.
For example, \cite{dellaert_approximate_2003} investigate uncapacitated lot sizing for a single item on a single machine and use a Markov decision model and dynamic programming to generate efficient lot sizing policies. Also \cite{Grubstrom_A_stochastic_2003} use dynamic programming, but aim to optimize the expected net present value for a multi-level capacitated lot sizing problem. Stochastic dual dynamic programming is used by \cite{quezada_combining_2022} in order to solve the uncapacitated version of the problem. On top of stochastic customer demand, the authors additionally investigate uncertain costs. \cite{van_hezewijk_using_2022} apply proximal policy optimization to tackle a lot sizing problem including multiple items, capacity restrictions and uncertain stationary demand.

Sequential decision problems under uncertainty can also be tackled by solving stochastic mixed integer programming models. For example, \cite{tunc_extended_2018} investigate extended formulations for stochastic lot sizing and present a novel cut generation approach. Demand scenarios are used by \cite{thevenin_material_2021}, who propose a two-stage, as well as a multi-stage stochastic programming model for multi-item multi-echelon capacitated lot sizing. Customer demand is considered as uncertain. A heuristic solution approach for the multi-stage decision problem is proposed and the resulting production plans are evaluated by means of a rolling-horizon scheme. As the authors point out, it is extremely challenging to directly solve the multi-stage model, due to its exponential growth in the number of decision stages. Making use of decomposition methods, such as progressive hedging, can help to resolve this issue.

Progressive hedging was initially proposed by \cite{rockafellar_scenarios_1991} and is a decomposition method for multi-stage stochastic programs iteratively solving the scenario subproblems and estimating the dual variables associated with the non-anticipativity constraints. They initially called the methodology “Scenario and Policy aggregation", as they aggregate information contained in the scenarios at various levels. \cite{wallace_structural_1991} further analyze the approach and provide some structural properties of the method.
Progressive hedging has been effectively utilized in several areas where stochasticity is important, such as network design, portfolio optimization, energy markets and production planning. \cite{jornsten_scenario_1994} apply progressive hedging to single resource production planning under demand uncertainty. \cite{hvattum_using_2008} use the progressive hedging algorithm to solve the stochastic inventory routing problem with scenario trees. \cite{kazemi_zanjani_scenario_2013} develop a scenario decomposition approach based on progressive hedging to solve a real-world production planning problem arising in sawmills. They consider randomness in demand, as well as in yield. \cite{zhao_multi-objective_2019} investigate pharmaceutical production planning for clinical studies in a multi-objective setting, minimizing production time and operational costs. They apply progressive hedging to tackle a multi-stage stochastic model and perform a numerical study determining optimal production quantities of clinical drugs. \cite{peng_progressive_2019} apply progressive hedging to the integrated production planning and scheduling problem in a multi-stage setting. They provide computational results when comparing their algorithm to commercial solvers and demonstrate the advantage of multi-stage modelling over two-stage and deterministic approaches. 

The progressive hedging concept can also be utilized to develop efficient heuristics. \cite{haugen_progressive_2001} cast the progressive hedging algorithm in a metaheuristic framework for single-item single-echelon lot sizing in a multi-stage stochastic setting. They solve the resulting scenario subproblems heuristically. Also \cite{crainic_progressive_2011} propose a metaheuristic based on the progressive hedging approach and apply it to stochastic network design. \cite{alvarez_inventory_2021} investigate the stochastic inventory routing problem with uncertain supply and demand. They present a two-stage stochastic optimization model and propose a progressive hedging inspired heuristic method to handle large instances. In \cite{sarayloo_integrated_2023} the progressive hedging approach is combined with a learning based metaheuristic in order to solve a capacitated network design problem under uncertain demand. \citet{kermani_progressive_2024} develop a progressive hedging-based matheuristic for solving the production routing problem under demand uncertainty.

Finally, \cite{thevenin_stochastic_2022} develop a hybrid method of progressive hedging and stochastic dual dynamic programming to solve multi-echelon lot sizing with component substitution in a multi-stage setting. In contrast to this work, setup carry-over decisions and external demand for components are not investigated. The latter influences the multi-stage model, because the non-anticipativity constraints for inventory and backlog variables are not implicitly given when component demands are considered. The additional decisions, whether components should be further processed to end items or be used to satisfy external demand, are represented by the corresponding inventory and backlog variables. Therefore, non-antipacitivity constraints are explicilty needed.

\subsection{Contributions}

The contributions of our work can be summarized as follows. 
\begin{itemize}
    \item First, we investigate the multi-stage stochastic decision process of multi-item multi-echelon lot sizing with setup carry-over, which represents a practically relevant extension for 24/7 production systems. We not only consider uncertain external demand for end items, but also for components along all BOM levels, motivated by industries that heavily rely on spare parts in the repair and maintenance process. 
    \item Second, we propose a progressive hedging approach to solve the multi-stage problem in the static-dynamic decision framework based on scenario trees. We investigate the crucial penalty parameter included in the scenario subproblems, which guides the subproblem solutions towards a consensus. There are no general rules available for the choice of this parameter and, therefore, our analysis provides valuable insights in the behavior of the progressive hedging approach for the multi-stage lot sizing problem at hand. 
    \item For the case of binary first-stage decision variables, the basic variant of PH is known to not necessarily converge to the global optimum. We therefore adapt metaheuristic adjustment strategies to guide the algorithm to a consensus more efficiently. The strategies significantly improve the convergence behaviour of PH and lead to higher quality solutions compared to the base version.
    \item Furthermore, we experiment with two different procedures to calculate the implementable solution at each iteration. First, we compute the average over the subproblem solutions, which is a common choice. Second, we also propose an algorithm variant making use of a majority voting instead of averaging. We demonstrate how the PH algorithm behaviour differs between the two variants.
    \item Finally, we compare the progressive hedging approach to a compact implicit model formulation and demonstrate the potential of progressive hedging. For several instances the tuned PH algorithm obtains solutions with a cost within 1\% of the  compact model solution cost in shorter runtimes. 
\end{itemize}  

\section{Problem formulation and mathematical model}
\label{sec:problem_formulation_and_model}

In Section \ref{sec:problem_formulation}, we define the problem setting for the multi-item multi-echelon capacitated lot sizing problem under demand uncertainty. Section \ref{sec:static-dynamic_decision_framework} then describes the static-dynamic decision framework for the multi-stage setting. A compact mathematical model for the proposed lot sizing problem is formulated in Section \ref{sec:mathematical_model}. Section \ref{sec:implicit_model} provides an implicit version of the compact model and Section \ref{sec:partial_implicit_model} describes how the complexity of the implicit model can be reduced by considering partial models including scenario path subsets. 

\subsection{Problem formulation}
\label{sec:problem_formulation}
In this work, we study a multi-stage stochastic decision problem, namely multi-item multi-echelon capacitated lot sizing with setup carry-over under demand uncertainty. The task comprises the production of multiple items, represented by the set $\mathcal{I}$, during the discrete time periods of the considered planning horizon $\mathcal{H}$, on the available resources $\mathcal{K}$. The aim is to efficiently meet the external and internal demand for items in order to minimize the overall cost. The notation used for the relevant problem parameters is presented in Table \ref{tab:modelparams}.
\begin{table}[!ht]
\begin{tabular*}{\hsize}{@{\extracolsep{\fill}}ll@{}}
\hline
Index set & Definition\\
\hline
$t, \tau \in \mathcal{H} = \{1,\ldots,T\}$ & Planning horizon \\
$i,j \in \mathcal{I} = \mathcal{I}_e \cup \mathcal{I}_c$ & Items \\
$\mathcal{I}_c = \{1,\ldots,N_c\}$ & Components \\
$\mathcal{I}_e = \{{N_c\!+\!1},\ldots,{N_c\!+\!N_e}\}$ & End items \\
$k \in \mathcal{K} = \{1,\ldots,J\}$ & Resources \\
$ \mathcal{I}_k$ for $k \in \mathcal{K}$ & Items that require resource $k$ \\
\hline
Parameter & Definition \\
\hline
    $D_{it}^{\phi}$  & Demand of item $i$ in period $t$ for scenario path $\phi$\\
    $R_{ij}$ & BOM - number of units of item $i$ required to \\ &produce one unit of item $j$ \\
    $L_i$ & Lead time of item $i$ \\
    $\hat{I}_{i0}$ & Initial inventory of item $i$ \\
    $C_{kt}$ & Capacity of resource $k$ in period $t$ \\
    $s_{i}$ & Setup cost of item $i$ \\
    $t_{i}$ & Setup time of item $i$ \\
    $p_{i}$ & Production time of item $i$ \\
    $v_i$ & Production cost of item $i$ \\
    $h_i$ & Holding cost of item $i$ \\
    $b_i$ & Backlog cost of end item $i$ \\
    $e_i$ & Lost Sale cost of end item $i$ (backlog at period $T$) \\
\hline 
\end{tabular*}
\caption{Model parameters}
\label{tab:modelparams}
\end{table}

In multi-echelon production systems, end items $\mathcal{I}_e$ are produced from several components $\mathcal{I}_c$ following a predefined BOM. The BOM is defined by a matrix $R_{ij}$ holding the information on the required number of units of component $i$ needed for the production of one unit of another component or end item $j$.  In each period of the discrete time horizon $\mathcal{H}$, external customer demand occurs for end items as well as for components, which is typical in industries that rely on the use of spare parts \citep{dodin_bombardier_2023}. Components therefore face both internal demand for the production of end items and external demand to satisfy customer orders. 

In order to describe the multi-stage decision process, we follow the notation of \cite{thevenin_stochastic_2022}. Customer demand in each period $t$ is stochastic and represented by discrete scenarios $\omega_t \in \Omega_t$. A set of scenarios $\omega_t$ for each $t \in \mathcal{H}$ is called a scenario path and is denoted by $\phi = (\omega_1, \ldots, \omega_T)$. A scenario path $\phi$ represents the realization of the uncertain demand along the planning horizon $\mathcal{H}$ and has probability $\sigma_{\phi}$. $D_{it}^{\phi}$ denotes the demand for item $i$ in period $t$ for scenario path $\phi$. Partial scenario paths are denoted by $\phi_{[t]}$ and contain demand realizations of $\phi$ up to period $t$. Let $\mathcal{D}_{\mathcal{I}[t]}^{\phi_{[t]}}$ denote the matrix of demand realizations $D_{it}^{\phi}$ corresponding to scenario path $\phi$ for all items $i \in \mathcal{I}$ up to period $t$.

We consider the static-dynamic decision framework, where setup decisions are fixed in the initial period, but setup carry-overs and production quantities, as well as component usage represented by inventory and backlog decisions can be adjusted in each period $t$, based on the demand information that is already available in the partial scenario path $\phi_{[t]}$. In Section \ref{sec:static-dynamic_decision_framework}, the scenario tree representation and the static-dynamic decision framework are further explained. The customer demand in a specific time period can be met by items in the inventory or the finished items of the respective period. Initially, there is an inventory of $\hat{I}_{i0}$ units available for item $i$. 

In this work, each item $i$ is assigned to exactly one resource $k$ out of the available resources $\mathcal{K}$. However, several items require the same resource. The set of items that are produced with resource $k$ is denoted by $\mathcal{I}_k$
. Before production can start, resources need to be set up. Setting up item $i$ on its necessary resource incurs a setup cost of $s_i$ and requires $t_i$ units of setup time. However, one setup per resource can be carried over to the next period and therefore remains intact. Producing one unit of item $i$ incurs a production cost of $v_i$. On top of that, an extra $p_i$ production time units are used on the necessary resource. For a time period $t$, $C_{kt}$ units of capacity can be used on resource $k$.

An item $i$ has a lead time $L_i$ before it can be further processed or used to satisfy external demand. Lead times are assumed to be smaller than or equal to 1. Longer lead times are rarely seen in big bucket lot sizing models and would unnecessarily complicate the capacity constraints.

Holding one unit of item $i$ in inventory for one period incurs a holding cost of $h_i$. If the available amount fails to meet the demand for item $i$, backlog cost of $b_i$ per unit and period is incurred until the demand is met. Any unsatisfied backlog for item $i$ in the final period of the planning horizon is treated as lost sales, which cannot be fulfilled later and carries a cost of $e_i$ per missing unit, which usually exceeds the backlog cost.

Decisions must be made on the production setups, which will be fixed over the whole planning horizon, as well as the setup carry-over, production quantity, inventory and backlog variables respecting non-anticipativity for all items in the face of uncertain demand.

\subsection{Static-dynamic decision framework}
\label{sec:static-dynamic_decision_framework}

\begin{figure}[ht]
    \centering  \includegraphics[width=\linewidth]{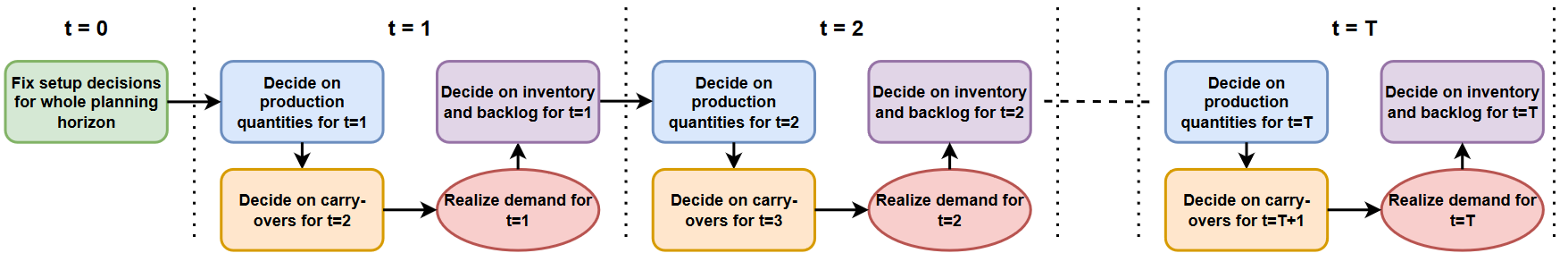}
    \caption{Static-dynamic decision framework}
    \label{fig:static-dynamic_decision_framework}
\end{figure}

We consider the static-dynamic decision framework. Under this framework, setup decisions are fixed for the whole planning horizon in the initial decision stage, while setup carry-overs, production quantities, inventories and backlogs are dynamically determined for each period, based on the most recent demand information. While production quantities for period $t$ are determined before the realization of the external demand in period $t$, the inventory and backlog decisions for period $t$ are made afterwards. Setup carry-over decisions are made simultaneously with the quantity decisions, since they define a partial schedule, which should be fixed as soon as production starts. While inventory and backlog for end items follow directly from the realized demand, additional decision-making is required, when also components face external demand. The decisions on inventory and backlog for components implicitly specify whether the respective external demand is immediately satisfied in order to reduce backlog or if components are stored in inventory in order to further turn them into end items in later periods. Even though instant fulfillment of external demand for components temporarily reduces backlog costs, the parts might be missing for end item production in subsequent periods and thus incur even higher backlog costs. Figure \ref{fig:static-dynamic_decision_framework} displays the underlying sequential decision process. 

\begin{figure}[ht]
    \centering
    \includegraphics[width=0.6\linewidth]{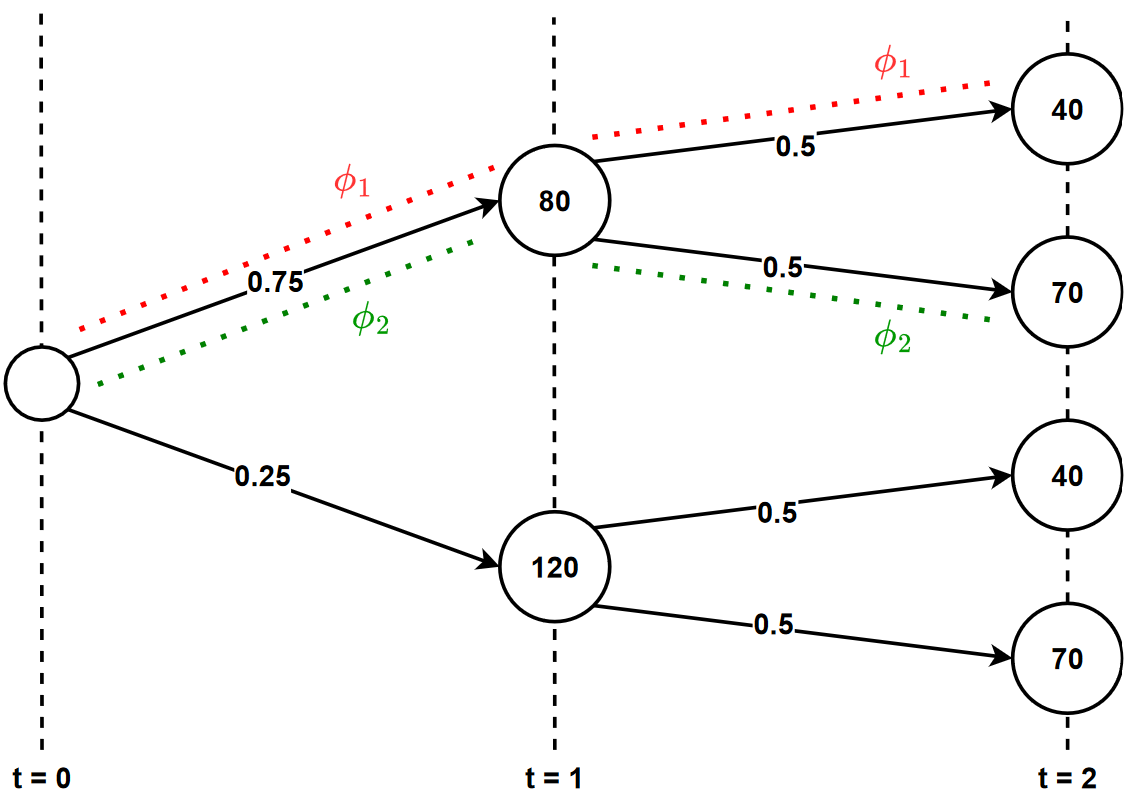}
    \caption{Example scenario tree for two periods}
    \label{fig:scenario_tree}
\end{figure}

When stochastic demand is represented by discrete scenarios in a multi-stage setting, so-called scenario trees are required. The scenario tree is the collection of all considered scenario paths $\phi \in \Phi$ and graphically demonstrates that several paths $\phi$ share a partial path $\phi_{[t]}$ from the beginning of the planning horizon up to period $t$, along which they cannot be distinguished. Figure \ref{fig:scenario_tree} shows an example scenario tree for two time periods. The scenario paths $\phi_1$ and $\phi_2$ are indistinguishable until time period $t=1$, since the partial scenario paths $\phi_{1[1]}$ and $\phi_{2[1]}$ are identical. Decisions related to these scenario paths need to be identical up to period $t=1$, because they rely on the same available information. Future realizations of the respective scenario paths cannot be anticipated and therefore may not influence the decisions -- they are subject to the so-called non-anticipativity constraints. We recall that $\mathcal{D}^{\phi_{[t]}}_{\mathcal{I}[t]}$ contains the matrix of all demand realizations $D_{it}^{\phi}$ corresponding to scenario path $\phi$ for each item $i \in \mathcal{I}$ and period $t \in \mathcal{H}$. For each period $t \in \mathcal{H}$ and scenario path $\phi \in \Phi$ we define the set $\mathcal{N}(\mathcal{D}^{\phi_{[t]}}_{\mathcal{I}[t]})$ that contains all scenario paths that are indistinguishable from $\phi$ up to period $t$.

\subsection{Mathematical model}
\label{sec:mathematical_model}

In the following we formulate the multi-stage stochastic programming model based on the models presented in \cite{thevenin_stochastic_2022} and \cite{SCHLENKRICH2024109379}. Table \ref{tab:modelvariables} presents the decision variables.
\begin{table}
\begin{tabular*}{\hsize}{@{\extracolsep{\fill}}ll@{}}
    \hline
    Variable & Definition \\
    \hline
    $Y_{it}$ & Binary variable taking value 1 iff there is a setup of item $i$ \\&in period $t$ \\
    $Z_{it}^{\phi}$ & Binary variable taking value 1 iff there is a setup carry-over of \\& item $i$ from period $t-1$ to period $t$ for scenario path $\phi$\\
    $Q_{it}^{\phi}$ & Production quantity of item $i$ in period $t$ for scenario path $\phi$\\
    $I_{it}^{\phi}$ & Inventory of item $i$ in period $t$ for scenario path $\phi$\\
    $B_{it}^{\phi}$ & Backlog of item $i$ in period $t$ for scenario path $\phi$ \\
    \hline
    \end{tabular*}
    \caption{Decision variables}
    \label{tab:modelvariables} 
\end{table}

\begin{equation}
\label{objective_function}
    \min \sum_{i \in \mathcal{I}_e \cup \mathcal{I}_c} \sum_{t \in \mathcal{H}} s_i Y_{it} +
    \sum_{\phi \in \Phi} \sigma_{\phi} \left (\sum_{i \in \mathcal{I}_e \cup \mathcal{I}_c} \left(\sum_{t \in \mathcal{H}} (v_i Q_{it}^{\phi} + h_i I_{it}^{\phi}) + \sum_{t=1}^{t=T-1} b_i B_{it}^{\phi} + e_i B_{iT}^{\phi} \right)\right )
\end{equation}
subject to
\begin{align}
\label{end_item_balance}
    \sum_{\tau=1}^{t-L_i}Q_{i\tau}^{\phi} + I_{i0}^{ \phi} - \sum_{\tau=1}^t D_{i\tau}^{\phi} - I_{it}^{\phi} + B_{it}^{ \phi} &= 0 &&\forall i \in \mathcal{I}_e, t \in \mathcal{H}, {\phi \in \Phi} \\
    \label{component_balance}
    \sum_{\tau=1}^{t-L_i}Q_{i\tau}^{\phi} + \hat{I}_{i0} - \sum_{\tau=1}^t(\sum_{j \in \mathcal{I}_e \cup \mathcal{I}_c} R_{ij} Q_{j\tau}^{\phi}) - \notag \\
    I_{it}^{ \phi} + B_{it}^{ \phi} - \sum_{\tau=1}^t D_{i\tau}^{\phi} &= 0 &&\forall i \in \mathcal{I}_c, t \in \mathcal{H}, { \phi \in \Phi}\\
\label{component_backlog}
B_{it}^{\phi} &\leq \sum_{\tau = 1}^t D_{i\tau}^{\phi}&&\forall i \in \mathcal{I}_c, t \in \mathcal{H}, { \phi \in \Phi}\\
\label{setup_or_carryover_constraint}
    Q_{it}^{\phi} - M_{it} (Y_{it} + Z_{it}^{\phi}) &\leq 0 &&\forall i \in \mathcal{I}, t \in \mathcal{H}, { \phi \in \Phi} \\
\label{one_carryover_per_resource}
    \sum_{i \in \mathcal{I}_k} Z_{it}^{\phi} &= 1 &&\forall t \geq 2, k \in \mathcal{K}, { \phi \in \Phi} \\
\label{carry_over_feasibility1}
    Z_{it}^{\phi} - Y_{i,t-1} - Z_{i,t-1}^{\phi} &\leq 0 && \forall i \in \mathcal{I}, t \geq 2, { \phi \in \Phi} \\
\label{carry_over_feasibility2}
    Z_{it}^{\phi} + Z_{i,t-1}^{\phi} - Y_{i,t-1} + Y_{j,t-1} &\leq 2 &&\forall i \in \mathcal{I}, j \neq i \in \mathcal{I}_k, t \geq 2, { \phi \in \Phi} \\
\label{capacity_contraint}
    \sum_{i \in \mathcal{I}_k}(t_i Y_{it} + p_i Q_{it}^{ \phi}) &\leq C_{kt} &&\forall t \in \mathcal{H}, k \in \mathcal{K}, { \phi \in \Phi} \\
\label{non-anticipativity_carry-over}
    Z_{i,t+1}^{\phi} &= Z_{i,t+1}^{\phi'} &&\forall i \in \mathcal{I}, t \in \mathcal{H}, \notag \\
    & &&\phi \in \Phi, \phi' \in \mathcal{N}(\mathcal{D}_{\mathcal{I}[t]}^{\phi_{[t]}})\\
\label{non-anticipativity_quantity}
    Q_{i,t+1}^{\phi} &= Q_{i,t+1}^{\phi'} &&\forall i \in \mathcal{I}, t \in \mathcal{H}, \notag \\
    & &&\phi \in \Phi, \phi' \in \mathcal{N}(\mathcal{D}_{\mathcal{I}[t]}^{\phi_{[t]}})\\
\label{non-anticipativity_inventory}
    I_{it}^{\phi} &= I_{it}^{\phi'} &&\forall i \in \mathcal{I}, t \in \mathcal{H}, \notag \\ & &&\phi \in \Phi, \phi' \in \mathcal{N}(\mathcal{D}_{\mathcal{I}[t]}^{\phi_{[t]}}) \\
\label{non-anticipativity_backlog}
    B_{it}^{\phi} &= B_{it}^{\phi'} &&\forall i \in \mathcal{I}, t \in \mathcal{H}, \notag \\ & &&\phi \in \Phi, \phi' \in \mathcal{N}(\mathcal{D}_{\mathcal{I}[t]}^{\phi_{[t]}})\\
\label{domain_setup}
    Y_{it} &\in \{0,1\} &&\forall i \in \mathcal{I}, t \in \mathcal{H} \\
\label{domain_carryover}
 Z_{it}^{\phi} &\in \{0,1\} &&\forall i \in \mathcal{I}, t \in \mathcal{H}, \phi \in \Phi \\
\label{domain_quantity}
Q_{it}^{ \phi} &\geq 0 &&\forall i \in \mathcal{I}, t \in \mathcal{H}, \phi \in \Phi \\
\label{domain_inventory_quantity}
I_{it}^{\phi} &\geq 0 &&\forall i \in \mathcal{I}, t \in \mathcal{H}, \phi \in \Phi \\
\label{domain_backlog}
B_{it}^{\phi} &\geq 0 &&\forall i \in \mathcal{I}, t \in \mathcal{H}, \phi \in \Phi.
\end{align}

The objective function \eqref{objective_function} aims to minimize the total cost over the considered planning horizon. The total cost include the setup costs, as well as the expected production, inventory, backlog and lost sales costs over all scenario paths $\phi \in \Phi$. Lost sales denote the unmet backlog of the last period of the planning horizon. Constraints \eqref{end_item_balance} and \eqref{component_balance} guarantee the balance of inventory and backlog for end items, as well as for components. This balance considers the required lead times and realizations of the uncertain customer demand. Since we consider external customer demands for components, constraints \eqref{component_backlog} make sure that the component backlog cannot exceed the respective external demand. This is necessary to comply with the assumption that the production of an item can only start when the required components are ready. In the absence of external component demand, the respective backlog variables are bound to zero, which simplifies the model. Constraints \eqref{setup_or_carryover_constraint} only allow production of an item $i$ in period $t$ if either a setup was performed or the setup for the respective item was carried over from period $t-1$. Equations \eqref{one_carryover_per_resource} make sure that on every resource only one set up state can be transferred to the next time period. The relation between setups and carry-overs is imposed by two sets of constraints. For the case of item $i \in \mathcal{I}$ being carried over on its needed machine $k \in \mathcal{K}_i$ in period $t \in \mathcal{H}$, Constraints \eqref{carry_over_feasibility1} guarantee that either a setup for item $i$ was performed in the previous period $t-1$ or the carry-over of the setup had already happened one period earlier, so from $t-2$ to $t-1$. Consider the latter case and assume the setup was already carried over from period $t-2$ to period $t-1$. Moreover, if a setup was made for another item $j \neq i \in \mathcal{I}$ on the required resource of item $i$ in period $t-1$, Constraints \eqref{carry_over_feasibility2} ensure that item $i$ can only be consecutively carried over to period $t$ if a new setup for item $i$ in period $t-1$ is performed. Constraints \eqref{capacity_contraint} limit the consumed resource capacity, caused by setup, as well as production activities. Since we consider the multi-stage setting, we must guarantee that decisions corresponding to scenario paths that are indistinguishable up to a period $t \in \mathcal{H}$ do not differ from each other until period $t$. These constraints \eqref{non-anticipativity_carry-over}-\eqref{non-anticipativity_backlog} are called non-anticipativity constraints because they make sure that knowledge that is not yet available cannot be used for decision making. Quantity decisions rely on the realized demand information of the previous period and therefore consistency needs to be guaranteed for period $t+1$. Inventory and backlog decisions are made immediately after demand is revealed and are therefore based on the demand of the same period $t$. Finally, \eqref{domain_setup}-\eqref{domain_backlog} specify the variable domains.

\subsubsection{Implicit model}
\label{sec:implicit_model}
The proposed multi-stage optimization model \eqref{objective_function}-\eqref{domain_backlog} contains many variables that are forced to have the same value by the non-anticipativity constraints. In order to make the model more efficient, variables that are supposed to have the same value can be grouped and replaced by one single variable. This not only reduces the number of variables but also removes the non-anticipativity constraints from the model. Following \cite{thevenin_material_2021} we replace the variables $\{Z_{i,t+1}^{\phi}| D_{it}^{\phi} = D(1,\ldots,t) \; \forall i \in \mathcal{I}, t \in \mathcal{H}, \phi \in \Phi\}$ by $Z_{i,t+1}^{D(1,\ldots,t)}$ and $\{Q_{i,t+1}^{\phi}| D_{it}^{\phi} = D(1,\ldots,t) \; \forall i \in \mathcal{I}, t \in \mathcal{H}, \phi \in \Phi\}$ by $Q_{i,t+1}^{D(1,\ldots,t)}$. Analogously we replace $\{I_{it}^{\phi}|D_{it}^{\phi} = D(1,\ldots,t) \;\forall i \in \mathcal{I}, t \in \mathcal{H}, \phi \in \Phi\}$ by $I_{it}^{D(1,\ldots,t)}$ and $\{B_{it}^{\phi}|D_{it}^{\phi} = D(1,\ldots,t) \;\forall i \in \mathcal{I}, t \in \mathcal{H}, \phi \in \Phi\}$ by $B_{it}^{D(1,\ldots,t)}$.

\subsubsection{Partial implicit model with scenario subset}
\label{sec:partial_implicit_model}
\begin{figure}
    \centering
    \includegraphics[width=0.8\linewidth]{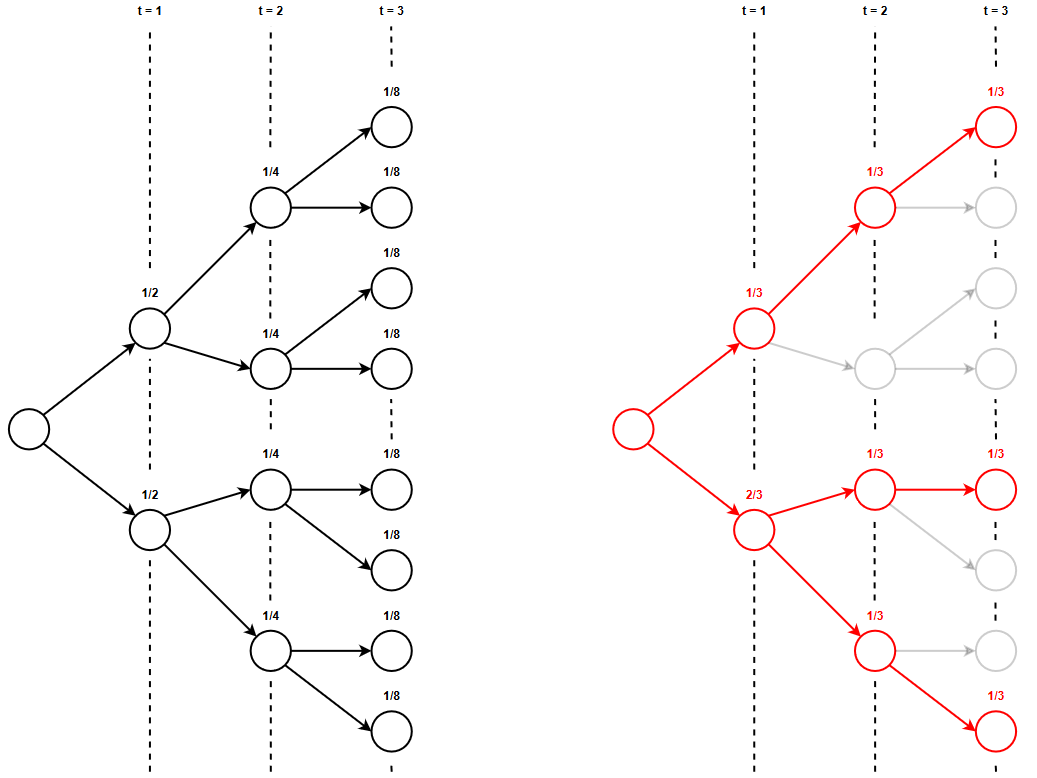}
    \caption{A partial scenario tree with adjusted scenario probabilities.}
    \label{fig:partial_implicit_model}
\end{figure}

The implicit model introduced in Section \ref{sec:implicit_model} is still a very challenging multi-stage stochastic optimization model. Even though it can be solved more efficiently than the compact model including all non-anticipativity constraints, solving the implicit model might still be hard for large instances including many scenario paths. A simple approach to reduce the complexity of the model is to decrease the number of included scenario paths. Instead of solving the model over all scenario paths $\phi \in \Phi$ we can define a subset $\Phi^{*} \subseteq \Phi$ and create a smaller version of the implicit model only including the scenario paths $\phi \in \Phi^{*}$. The first stage solutions to the simplified implicit model are also feasible for the full model including all scenario paths $\phi \in \Phi$. Therefore, solving the simplified model can be seen as a heuristic approach to solve the full formulation. Note that the probabilities for the selected scenario paths $\phi \in \Phi^{*}$ need to be adjusted, compared to the full scenario tree, as shown in Figure \ref{fig:partial_implicit_model}.

\section{Progressive hedging approach}
\label{sec:progressive_hedging}

Progressive hedging, initially proposed by \cite{rockafellar_scenarios_1991}, is a decomposition method for multi-stage stochastic programs. It iteratively solves the deterministic scenario subproblems by relaxing the non-anticipativity constraints, but applies a penalty to the objective function in case these constraints are violated. We follow the presentation of \cite{thevenin_stochastic_2022} but note that the non-aticipativity constraints for inventory and backlog variables are not implicitly given in our work, since we consider external demand for components. 

The multi-stage stochastic optimization model \eqref{objective_function}-\eqref{domain_backlog} can be decomposed by scenario path $\phi \in \Phi$. In order to do this, we equip the first stage decision variables $Y_{it}$ with scenario path indices, i.e., $Y_{it}^{\phi}$. We then add non-anticipativity constraints that force them to be equal for all scenario paths, i.e., $Y_{it}^{\phi} = Y_{it}^{\phi'}$ for all $i \in \mathcal{I}, t \in \mathcal{H}, \phi, \phi' \in \Phi$. By relaxing the non-anticipativity constraints \eqref{non-anticipativity_carry-over}-\eqref{non-anticipativity_backlog}, as well as the newly introduced constraints for the first stage variables the problem decomposes into a deterministic lot sizing problem for each scenario path $\phi \in \Phi$, denoted as scenario subproblem $\mathcal{P}_{\phi}$.

A solution for $\mathcal{P}_{\phi}$ is denoted as the vector $x_{\phi}$ that contains the optimal values for the decision variables, denoted by the set $\nu = \{Y_{it}^{\phi},Z_{it}^{\phi},Q_{it}^{\phi},I_{it}^{\phi},B_{it}^{\phi} \; \forall i \in \mathcal{I}, t \in \mathcal{H}\}$. Determining the subproblem solutions $x_\phi$ for all scenario paths $\phi \in \Phi$ might lead to a non-implementable solution. This is because there is no guarantee for decisions at period $t \in \mathcal{H}$ to be equal for two scenario paths $\phi$ and $\phi'$ that are indistinguishable up to time $t$. Even though the decisions are based on the same available information at time $t$, the decisions might differ, which is not feasible in the static-dynamic setting. In order to solve this issue, implementable decisions, also called consensus decisions, are calculated at each stage based on the subproblem solutions. A consensus can be calculated in several ways. However, the most common approach is to take the average over all respective subproblem solutions. Later, we present a version of the algorithm that uses a majority voting instead. Deviations from these implementable decisions are then penalized in the objective functions of the deterministic scenario subproblems. The magnitude of the penalty has an important impact on the convergence behaviour of the algorithm and is mainly influenced by the penalty parameter $\rho$. Therefore, we will investigate several choices for $\rho$. The scenario subproblems are then solved iteratively until the subproblem solutions converge to the implementable solution. Several stopping criteria can be defined. However, the most common one is to terminate the PH algorithm, as soon as the first stage variables, i.e., the setup decisions, have converged to the consensus. Algorithm \ref{alg:progressive_hedging} presents the procedure in more detail.

\begin{algorithm}
\caption{Progressive hedging algorithm}\label{alg:progressive_hedging}
\begin{algorithmic}[1]

\State $k \gets 1$
\State $\Lambda_{x}^k \gets 0$
\State $x_{\phi}^k \gets$ Solve scenario subproblem $\mathcal{P}^{k}_{\phi}$ for all scenario paths ${\phi} \in \Phi$
\State $\Tilde{x} \gets$ Calculate implementable solution based on solutions for scenario subproblems $x^k_{\phi}$ and $\mathcal{N}(\mathcal{D}^{\phi_{[t-1]}}_{\mathcal{I}[t-1]})$
\State $\Lambda_{x}^k \gets \Lambda_{x}^{k-1} + \rho(x^k - \Tilde{x}) $ Update Lagrange multipliers 
\State Update penalty term in objective function of scenario subproblems $\mathcal{P}^{k}_{\phi}$ based on new implementable solution $\Tilde{x}$ and new Lagrangian multipliers $\Lambda_{x}^k$: $\sum_{x \in \nu} \Lambda_{x}^k(x-\Tilde{x}) + \frac{\rho}{2}(x-\Tilde{x})^2$
\If{$x^k$ converged to $\Tilde{x}$, meaning $(x^k - \Tilde{x}) \leq \epsilon$}
    \State return $\Tilde{x}$
\Else 
    \State $k \gets k+1$
    \State go to 3
\EndIf

\end{algorithmic}
\end{algorithm}

\subsection{Choosing penalty parameter $\rho$}
The penalty parameter $\rho$ is an important parameter influencing the speed of convergence of the PH algorithm. It affects the penalty term in the objective function of the scenario subproblems $\mathcal{P}^{k}_{\phi}$ twice:
\begin{enumerate}
    \item Linear effect: the value of $\rho$ has a linear impact on the update of the Langrange multipliers $\Lambda_{x}^k \gets \Lambda_{x}^{k-1} + \rho(x^k - \Tilde{x}) $;
    \item Quadratic effect: $\frac{\rho}{2}$ is also the coefficient of the quadratic penalty  $\frac{\rho}{2}(x-\Tilde{x})^2$.
\end{enumerate}
As \cite{thevenin_stochastic_2022} point out, there is no general rule for the most suitable choice of $\rho$, since it strongly depends on the characteristics of the problem at hand. As \cite{watson_progressive_2011} highlight, the value of $\rho$ does not need to be equal for all decision variables, but can be chosen individually based on the respective objective coefficients. This can be advantageous, because if the magnitude of $\rho$ is not close to the magnitude of the objective coefficient, this will result in a small change in the penalty term in the objective function between iterations \citep{peng_progressive_2019}. For this reason \cite{watson_progressive_2011} propose different strategies for choosing $\rho$, among which is one they call the \textit{cost proportional approach}. In this strategy the penalty parameter $\rho_x$ for decision variable $x$ is chosen as a $\lambda$-multiple of the respective objective function coefficent $\text{coef}_x$. We adapt this strategy in our work and evaluate different values for $\lambda$ in the set $\{0.1, 1, 10, 100\}$.

\subsection{Metaheuristic adjustment strategies}
\label{sec:metaheurisitc_adjustment_strategies}

For the PH version proposed in \cite{rockafellar_scenarios_1991}, and respectively the algorithmic procedure of Algorithm \ref{alg:progressive_hedging}, the authors show that the approach converges to a globally optimal solution for continuous problems. However, this is not implied when binary decision variables are considered. In order to guide the algorithm more efficiently towards a consensus we implement two adaptation strategies that are based on the strategies proposed by \cite{kermani_progressive_2024} and \cite{crainic_progressive_2011}. In the work of \cite{crainic_progressive_2011}, it was shown that heuristically adjusting penalty and cost parameters can lead the first-stage variables to a consensus more efficiently, especially for binary first-stage variables. The idea of the adjustment strategies applied in this work is twofold. 

First, the so-called global adjustment strategy should make a deviation from the overall consensus in the subproblems less likely in case the implementable solution has surpassed a certain threshold. If the implementable solution value of a setup variable $\bar{y}_{it}$ is larger than a threshold $\theta_H$ (with $0.5 < \theta_H < 1$), this means that a majority of the scenario subproblems agree to perform a setup for item $i$ in period $t$. On the contrary, if the implementable solution value of a setup variable $\bar{y}_{it}$ is lower than a threshold $\theta_L$ (with $0 < \theta_L < 0.5$), this means that a majority of the subproblems agree to not make a setup for the respective item in period $t$. Of course, during the PH procedure the fractional value of the binary decision variables can significantly fluctuate and also fall below a certain threshold after surpassing it in earlier iterations. However, to improve the convergence speed of the approach, parameters can be adapted based on these thresholds in order to push the scenario subproblem solutions in the direction of the consensus. \cite{kermani_progressive_2024} increase setup costs along all subproblems for an item $i$ in period $t$ in case the implementable solution for the respective setup falls below $\theta_L$ in order to make the appearance of a setup in a subproblem solution even less likely. Equivalently, they reduce setup costs for an item $i$ in period $t$ along all subproblems in case the implementable solution surpasses $\theta_H$ in order to support the creation of a setup in the last remaining subproblems that had no setup so far.

Second, the so-called local adjustment strategy modifies parameters for specific scenario subproblems, rather than adjusting them along all of them. This strategy is triggered when the absolute deviation of a scenario path solution significantly varies from the implementable solution. In order to measure this, a third threshold $\gamma_F$ is introduced, with $ 0.5 < \gamma_f < 1$. \cite{kermani_progressive_2024} then modify the setup costs for item $i$ in period $t$ of a single scenario subproblem similar to the global adjustment strategy, based on the respective scenario subproblem solution value. This should force the setup value to zero, in case the consensus is already zero among a large majority, and vice versa to one in case a large majority agrees to make a setup.

We note that the problem considered by \cite{kermani_progressive_2024} only has one item and no setup carry-overs, which makes the setup costs the main driver for the setup decisions. Since the problem setting in this work differs in several respects, instead of using the setup costs as the only tool to push setup decisions we also modify the penalty parameter $\rho$, which explicitly has the purpose to punish deviations from the implementable solution. This strategy follows the approach of \cite{crainic_progressive_2011}. As the authors note, when considering integer decision variables, PH does not necessarily converge to the global optimum. Since we tackle the binary case, the approach represents a heuristic and the multiplier modifications give an incentive to make or remove a setup when its value significantly differs from the current consensus. This incentive is gradually growing with an increasing penalty factor. In the following, we report the slightly modified adjustment strategies used in this work, which use both ideas -- modifying objective function cost parameters, as well as penalty parameters -- in order to guide the algorithm towards a consensus more efficiently:

\begin{enumerate}
    \item Global adjustment strategy: adjust the setup cost of a specific item $i$ and time period $t$, if the implementable decision variable $\bar{y}_{it}$ is fractional, depending on a predefined threshold with $0 < \theta_L < 0.5$ ($0.4$ is chosen) and $0.5 < \theta_H < 1$ ($0.6$ is chosen), where $\lambda_G > 1$ is the global modification rate:
    \begin{equation}
    s_{it}^k = 
    \begin{cases} 
        \lambda_G s_{it}^{k-1} & \text{if } \bar{y}_{it} < \theta_L \\
        \frac{1}{\lambda_G} s_{it}^{k-1} & \text{if } \bar{y}_{it} > \theta_H \\
        s_{it}^{k-1} & \text{otherwise.}
    \end{cases}
    \end{equation}
    
    \item Local adjustment strategy: modify the penalty parameter $\rho$ for the setup variable of an item in a certain period for a scenario subproblem $\phi$, if the distance to the implementable solution is above a threshold $\gamma_F$ with $0.5 < \gamma_F < 1$ ($0.8$ is chosen). The local modification rate is denoted by $\lambda_L$ (with $\lambda_L > 1$):
    \begin{equation}
    \rho_{it}^k(\phi) = 
    \begin{cases} 
        \lambda_L \rho_{it}^{k-1}(\phi) & \text{if } |y_{it}^{k-1} - \bar{y}_{it}| \geq \gamma_F \\
        \rho_{it}^{k-1}(\phi) & \text{otherwise.}
    \end{cases}
    \end{equation}


\end{enumerate}

\subsection{Majority voting for implementable solution}
\label{sec:majority_voting}

There are several ways to compute an implementable (consensus) solution based on the different subproblem solutions. Instead of averaging the respective variable values, a majority voting can be applied. This can be especially useful for the first-stage binary decision variables. In case the averaged variable lies within the interval $[0,0.5]$ the implementable solution for this variable is rounded to $0$, while it is rounded to $1$ if the value lies between $(0.5,1]$. 

This implementation changes the behaviour of the PH algorithm in the following sense. Assume that there are only three scenario paths (denoted by 1, 2 and 3) included in the PH algorithm. Consider the setup decision for item $i$ in period $t$ and the respective scenario subproblem solutions for this particular setup: $y_{it}^1, y_{it}^2, y_{it}^3$. Assume that scenario subproblems 1 and 2 decide to perform a setup, hence $y_{it}^1 = 1$ and $y_{it}^2 = 1$, but in scenario subproblem 3 no setup is performed, i.e., $y_{it}^3 = 0$. When using the averaging consensus calculation the implementable solution for item $i$ in period $t$ would then be $\bar{y}_{it} = 0.66$. Deviating from this consensus will be penalized in the next iteration of the PH algorithm. However, since the setup decisions are binary, all scenario subproblems will deviate from this implementable solution -- some more and some less. Performing a setup, and hence setting the variable to 1 will lead to a smaller violation (0.33) than not performing a setup, which leads to a violation of 0.66. 

Majority voting tackles this in a different way. For the case described above, the implementable solution is set to 1, since the majority of the scenario subproblems agree to perform a setup, i.e., $\bar{y}_{it} = 1$. This implies that in the following PH iteration scenario subproblems that perform a setup for item $i$ in period $t$ will not be penalized. On the other hand the penalty for the scenario subproblems that still do not perform a setup will be even larger.

\section{Computational study}
\label{sec:computational_study}
In this section, we perform a computational study of the proposed PH algorithm and the implicit multi-stage stochastic optimization model on well known multi-item multi-echelon benchmark instances. All models, as well as the PH algorithm, were implemented in Python using the CPLEX Python API and CPLEX 20.1.0.0 as MIP solver. All experiments were carried out on a  Quad-core X5570 Xeon CPU @2.93GHz with 48GB RAM. 

\subsection{Test instances}
\label{sec:test_instances}
\begin{figure}[!ht]
    \centering
    \includegraphics[width=10cm]{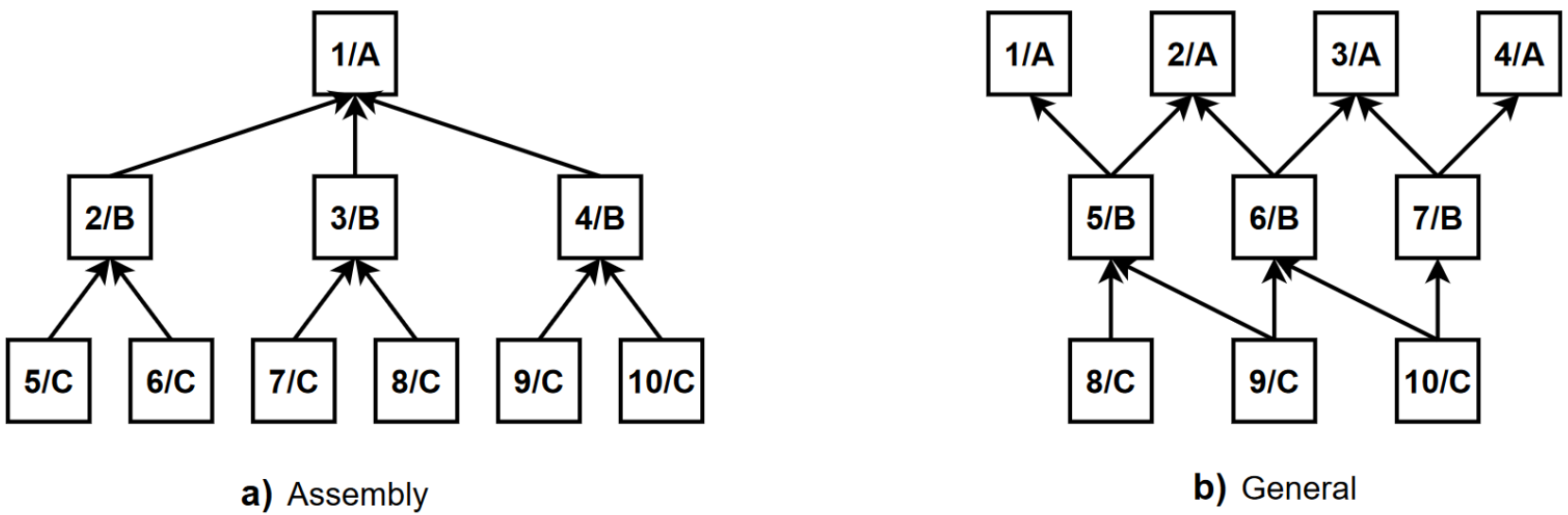}
    \caption{BOM for assembly (a) and general (b) production structures including resource assignments ($i / k$ - item $i$ is assigned to resource $k$)}
    \label{fig:BOM_structures}
\end{figure}

For our computational study we use the well known benchmark instances initially proposed by \cite{tempelmeier_derstroff_1996} as a basis and adjust the instances to our problem setting according to \cite{thevenin_stochastic_2022} and \cite{sereshti_managing_2024}. While we follow \cite{thevenin_stochastic_2022} regarding the generation of scenario trees, we use the procedure of \cite{sereshti_managing_2024} to generate external demands for components. The test instances comprise 7 time periods ($|\mathcal{H}| = 7$), 10 items ($|\mathcal{I}|=10$) and two BOM structures, namely an assembly and a general BOM structure. Figure \ref{fig:BOM_structures} shows the respective BOM structures, as well as the resource assignments of items at all levels. We consider instances including external demand only for end items, as well as instances including component demands. Moreover, we consider resource utilizations of 50\% and 90\% and three types of setup times, namely no setup times, short and long. The considered times are following test set B+ in \cite{tempelmeier_derstroff_1996}, where we take profile 1 and profile 3 for short and long times respectively. We consider backlog to holding cost ratios of 2 and 4 and further consider two variants for the holding costs, following \cite{thevenin_material_2021}, where the authors propose one variant with constant echelon holding costs and one with high holding costs for later echelons. In total this results in 96 test instances. 

As demand distribution we choose a modified version of the Lumpy demand distribution used in \cite{thevenin_material_2021}. We aim to model fluctuating demands with three levels, namely regular demand, zero demand and extremely high demand, extending the two levels of the Lumpy distribution in \cite{thevenin_material_2021}. Preliminary results have shown that this type of demand uncertainty has a high impact on the resulting setup patterns and therefore makes the problem setting especially challenging. The Lumpy demand distribution used in this work is defined as follows. With a probability of $\frac{1}{2}$ the demand is regular, being modeled by drawing from a Poisson distribution with mean $\frac{2}{3} F_{it} $, where  $F_{it}$ is the average demand provided in \cite{tempelmeier_derstroff_1996}. With a probability of $\frac{1}{3}$ there is no demand at all for an item in a certain period. Finally, with a probability of $\frac{1}{6}$ the demand is extremely high, being modeled by drawing from a Poisson distribution with mean $4 F_{it} $.

\subsection{Full implicit model results}
As a baseline we run the implicit multi-stage optimization model presented in Section \ref{sec:mathematical_model} on the 96 test instances described in Section \ref{sec:test_instances}. Table \ref{tab:Lumpy_full_implicit_model} reports the average runtimes in seconds and average obtained gaps in percent, as well as the percentage of instances, for which the implicit model either reached the time limit of three hours (Lim) or the solver ran out of memory (Mem). Runtimes and gaps are averaged over the instances, where the solver did not run out of memory. The obtained results show that the full implicit multi-stage model is already very challenging for three branches per scenario leading to a total of 2,187 scenario paths. In the case of 90\% utilization we are not able to find the optimal solution within the time limit of three hours in 33.3\% of the cases for end item demand only and in 62.5\% of the cases when components face external demand. In particular, the cases including component demands seem to be very challenging. Only for the case of $|\Omega| = 2$ we are able to obtain the optimal solution for all tested instances. For $|\Omega| = 4$ there are significant gaps after a runtime of 3 hours for at least half of the tested instances. In the case of external component demand and $|\Omega| = 4$ the time limit is reached for more than 90\% of the instances. In the case of $|\Omega| = 5$ the solver runs out of memory for at least 25\% of the instances and reaches the time limit for almost all remaining runs. This demonstrates the limitations of solving the multi-stage decision problem by means of a compact optimization model. Note that even though the considered values for the size of $\Omega$ represent a small sample for drawing from a random distribution, the exponential growth of the scenario tree makes the problem intractable. We therefore experiment with the proposed PH algorithm in the next section. 
\begin{table}
    \centering
    \resizebox{\textwidth}{!}{
    \begin{tabular}{c|c|r|r|r|r|r|r|r}
    \hline
    \multicolumn{3}{c}{}& \multicolumn{3}{|c}{End item demand only} & \multicolumn{3}{|c}{Component demand} \\
    \hline 
    Util. & $|\Omega|$ & $|\Phi|$ & Runtime (s) & Gap (\%) & Lim/Mem (\%)& Runtime (s) & Gap (\%) & Lim/Mem (\%)\\
    \hline
        \multirow{5}{*}{50\%} 
         & 2 & 128 & 67.5 & 0.00 & 0.0/0.0 & 92.5 & 0.00 & 0.0/0.0 \\
         & 3 & 2,187 & 2770.6 & 0.51 & 16.6/0.0 & 4274.6 & 0.02 & 4.1/0.0\\
         & 4 & 16,384 & 6504.9 & 10.24 & 50.0/0.0 & 10150.0 & 24.11 & 91.6/0.0\\
         & 5 & 78,125 & 8879.7 & 29.28 & 58.3/25.0 & 10728.0 & 56.68 & 66.6/29.1\\
         \hline
         \multirow{5}{*}{90\%} 
         & 2 & 128 & 1038.8 & 0.00 & 0.0/0.0 & 841.1 & 0.00 & 0.0/0.0\\
         & 3 & 2,187 & 4225.4 & 3.07 & 33.3/0.0 & 7628.6 & 1.80 & 62.5/0.0\\
         & 4 & 16,384 & 7377.2 & 19.49 & 50.0/0.0 & 10033.9 & 40.12 & 91.6/0.0\\
         & 5 & 78,125 & 8929.5 & 72.41 &  37.5/54.2 & 10800.0 & 66.24 & 62.5/37.5\\
         \hline
    \end{tabular}
    }
    \caption{Averaged runtimes and gaps for the full implicit multi-stage optimization model on 96 instances based on the \cite{tempelmeier_derstroff_1996} instances ($|\mathcal{I}| = 10, |\mathcal{H}| = 7$) for different numbers of scenarios per branch with a time limit of 3 hours (10800 s). Results are grouped by resource utilization and demand type and we average runtimes and gaps over the instances where the solver did not run out of memory.}
    \label{tab:Lumpy_full_implicit_model}
\end{table}

\subsection{Evaluating penalty parameter $\rho$ for PH}
In order to be able to put the performance of the PH algorithm into perspective, we run the following experiments for $|\Omega| = 2$, since it is the only case where we could obtain the optimal solution of the implicit model for all instances. First of all, we evaluate the penalty parameter $\rho$ for PH, which is the major driver for penalizing a deviation from the implementable solution and has a significant impact on the performance of the algorithm. More precisely, we investigate the $\lambda$-multiples of the objective function coefficient for each decision variable, because we follow the cost proportional approach of \cite{watson_progressive_2011}. We test values in $\{0.1, 1, 10, 100\}$ for $\lambda$. Let $coef_x$ denote the cost coefficient of a decision variable $x$, we then  calculate the respective penalty parameter as $\rho_x = \lambda coef_x$. 

Table \ref{tab:evaluation_rho_AVG_no_adjustments} shows the results for the evaluation of the $\lambda$-multiple for the PH algorithm. It becomes clear that for small $\lambda$-values the PH struggles to converge to a consensus within the time limit of three hours for a majority of the test instances. For $\lambda=0.1$ the algorithm does not converge for any problem instance considering external demand for components. In the case of only end item demands, PH with $\lambda=0.1$ converges for only 8.3\% of the instances with low utilization, and 20.8\% for high utilization instances. Also for $\lambda=1$ the percentages of converged PH runs are rather low with values between 8.3\% and 29.2\%. For those converged runs, however, the obtained cost delta to the optimal cost is very low and lies between 0.09\% and 0.83\%. For high $\lambda$-values, such as 10 and 100, the PH algorithm converges in almost all cases for low utilization of 50\% and even for the hard case of 90\% utilization, PH converges in between 66.6\% and 79.2\% of all runs. Also, the averaged runtimes obtained by higher $\lambda$-values are much smaller compared to the PH runs with small $\lambda$-values. However, applying a $\lambda$ of 10 or 100 leads to higher quality deltas between 2.10\% and 5.35\% to the optimal costs. It appears that the algorithm would benefit from a smaller $\lambda$-value in terms of solution quality, which, however, comes at the cost of extremely long runtimes. In order to counteract this effect, we test the proposed adjustment strategies and their effect on the PH algorithm behavior.
\begin{table}[ht]
    \resizebox{\textwidth}{!}{
    \centering
    \begin{tabular}{c|c|r|r|r|r|r|r}
    \hline
    \multicolumn{2}{c}{} & \multicolumn{3}{|c}{End item demand only} & \multicolumn{3}{|c}{Component demand}\\
    \hline 
    Utilization & $\lambda$ & Converged (\%) & Runtime (s) & $\Delta$ (\%) & Converged (\%) & Runtime (s) & $\Delta$ (\%)\\
    \hline
        \multirow{4}{*}{50\%} 
         & 0.1  & 8.3 & 2662.9 & 0.00 & 0.0 & - & - \\
         & 1  & 29.1 & 3992.3 & 0.47 & 29.2 & 7816.1 & 0.56 \\
         & 10  & 91.6 & 225.1 & 3.41 & 100.0 & 322.1 & 3.26 \\
         & 100  & 100.0 & 232.4 & 5.35 & 100.0 & 208.8 & 4.62\\
         \hline
         \multirow{4}{*}{90\%} 
         & 0.1  & 20.8 & 2041.1 & 0.00 & 0.0 & - & - \\
         & 1 & 25.0 & 1138.7 & 0.83 & 8.3 & 4264.0 & 0.09 \\
         & 10  & 75.0 & 401.7 & 3.24 & 79.2 & 1390.2  & 2.10\\
         & 100  & 79.2 & 636.5 & 4.46 & 66.6 & 1262.3 & 3.53 \\
         \hline
    \end{tabular}
    }
    \caption{Averaging consensus calculation: Evaluation of the $\lambda$-multiple in terms of total costs and runtime on a scenario tree with $|\Omega| = 2$. Reporting percentage of converged PH runs and respective average runtimes and average cost deltas to optimal costs averaged over the converged runs for each individual $\lambda$-value. Results are reported for the 96 test instances grouped by utilization and demand type.}
    \label{tab:evaluation_rho_AVG_no_adjustments}
\end{table}

\subsection{Evaluating penalty parameter $\rho$ for PH with adjustment strategies}
Since metaheuristic adjustment strategies have already been shown to have the potential to guide the PH algorithm towards a consensus more efficiently, we test the proposed global and local adjustment strategies described in Section \ref{sec:metaheurisitc_adjustment_strategies}. For thresholds and modifiers, we choose the following values, based on the experiments of \cite{kermani_progressive_2024}: $\theta_L = 0.4, \theta_H = 0.6, \gamma_F = 0.8, \lambda_G = 1.1, \lambda_L = 1.5$. We again test the same set of values for the $\lambda$-multiplier, namely $\{0.1, 1, 10, 100\}$.

Table \ref{tab:evaluation_rho_AVG} shows the results for the evaluation of the $\lambda$-multiplier for PH with adjustment strategies. It reports average runtimes and cost deltas for the 96 test instances grouped by utilization and demand type. For each grouping we highlight the best obtained cost delta to the optimal cost in bold. First of all, the most notable difference to the experiments without applying the metaheuristic adjustment strategies is that PH naturally converged for almost all instances across all $\lambda$-multipliers within the time limit of three hours. The minimum percentage of naturally converged runs over the four groupings for each $\lambda$-value was 83.3\% for $\lambda=0.1$, 87.5\% for $\lambda=1$, 91.6\% for $\lambda=10$ and 95.8\% for $\lambda=100$. For the few remaining instances a cyclic pattern was recognized and disrupted by significantly increasing $\rho$ by factor 10 to force convergence.

Second, the experiments clearly suggest the choice of $\lambda = 1$ as the most suitable multiplier for the penalty parameter $\rho$. For the investigated cases the resulting averaged gap to the optimal solution is less than or close to 1\%. As assumed the PH runtime decreases with an increase in $\lambda$, however this also comes with a decrease in solution quality. The resulting gap for the largest tested $\lambda$-value 100 lies between 3.12\% and 5.60\%. Choosing $\lambda = 0.1$ leads to longer runtimes and larger gaps to the optimal solution, compared to $\lambda=1$. 

\begin{table}[ht]
    \centering
    \begin{tabular}{c|c|r|r|r|r}
    \hline
    \multicolumn{2}{c}{} & \multicolumn{2}{|c}{End item demand only} & \multicolumn{2}{|c}{Component demand} \\
    \hline 
    Utilization & $\lambda$ &  Runtime (s) & $\Delta$ (\%) &  Runtime (s) & $\Delta$ (\%) \\
    \hline
        \multirow{4}{*}{50\%} 
         & 0.1  & 1046.1 & 0.76  & 1931.3   & 3.72\\
         & 1  & 549.7  & \textbf{0.36}  & 779.5  &  \textbf{0.50} \\
         & 10   & 239.7  & 2.35  & 390.5 & 3.06 \\
         & 100 & 184.5  & 4.18  & 191.3 & 4.35 \\
         \hline
         \multirow{4}{*}{90\%} 
         & 0.1  & 1331.6 & 2.51  & 2991.5 & 0.87 \\
         & 1  & 791.7  & \textbf{1.57}  & 1612.6 & \textbf{0.75} \\
         & 10   & 458.7   & 3.77  & 861.0  & 1.83 \\
         & 100   & 372.8 & 5.60  & 767.6  & 3.12  \\
         \hline
    \end{tabular}
    \caption{Averaging consensus calculation: Evaluation of the $\lambda$-multiple in terms of average runtime and cost delta to optimal cost on a scenario tree with $|\Omega| = 2$. Results are reported for 96 test instances grouped by utilization and component type.}
    \label{tab:evaluation_rho_AVG}
\end{table}

\subsection{Majority voting for consensus calculation}
Instead of averaging, the consensus solution can also be determined by a majority voting, as explained in Section \ref{sec:majority_voting}. The results for the majority voting consensus calculation can be found in Table \ref{tab:evaluation_rho_MV}. Most notable for the majority voting consensus calculation, every single PH run naturally converged within the time limit of three hours. We again report average runtimes and cost deltas over the 96 test instances grouped by resource utilization and demand type. For each grouping we highlight the best obtained cost delta to the optimal cost in bold.

The runtimes of the PH algorithm significantly decreased when using the majority voting consensus calculation across all $\lambda$-values when compared to averaging. Again, also for the majority voting setting, $\lambda = 1$ appears to be the most suitable choice, since it provides the smallest delta to the optimal costs. However, the average cost deltas are larger for majority voting than for averaging. The decrease in runtime therefore also comes with a decrease in solution quality. When compared to the averaging consensus calculation presented in Table \ref{tab:evaluation_rho_AVG}, a different pattern becomes visible. When choosing the majority voting, a high $\lambda$-value does not imply a large increase in cost delta to the optimal cost in three out of four cases, as it was the case when averaging was used for the consensus calculation. Also in terms of PH runtimes the majority voting leads to a different pattern. Instead of decreasing runtimes for higher $\lambda$-values, they increase. 
The PH algorithm with majority voting appears not to suffer from high penalties in terms of quality. However, it also does not benefit from them in terms of runtime. For PH with averaging this was clearly the case. The choice for $\lambda = 1$, however, remains the most suitable, independent of the consensus calculation.

\begin{table}
    \centering
    \begin{tabular}{c|c|r|r|r|r}
    \hline
    \multicolumn{2}{c}{} & \multicolumn{2}{|c}{End item demand only} & \multicolumn{2}{|c}{Component demand} \\
    \hline 
    Utilization & $\lambda$ & Runtime (s) & $\Delta$ (\%) & Runtime (s) & $\Delta$ (\%)\\
    \hline
        \multirow{4}{*}{50\%} 
         & 0.1  & 809.2 & 5.22 & 1651.4  & 7.83  \\
         & 1  & 324.4 & \textbf{1.35}  & 533.4   & \textbf{2.69} \\
         & 10  & 260.7  & 2.40  & 354.4 & 3.03 \\
         & 100  & 288.9  & 2.53  & 377.3 & 3.06 \\
         \hline
         \multirow{4}{*}{90\%} 
         & 0.1  & 905.4 & 7.79  & 1950.4  & 5.31\\
         & 1  & 394.5 & \textbf{4.35}  & 1012.4  & \textbf{2.48}\\
         & 10  & 436.3 & 4.63  & 600.0  & 2.50\\
         & 100  & 329.2 & 5.91  & 768.8  & \textbf{2.48} \\
         \hline
    \end{tabular}
    \caption{Majority voting consensus calculation: Evaluation of the $\lambda$-multiple in terms of average runtime and cost delta to optimal cost on a scenario tree with $|\Omega| = 2$. Results are reported for 96 test instances grouped by utilization and component type.}
    \label{tab:evaluation_rho_MV}
\end{table}

\subsection{Partial scenario tree evaluations}
In order to also address scenario trees with a larger number of scenarios $|\Omega|$ per tree node, which cannot be solved to optimality with the full implicit model for a majority of the instances, we experiment with the partial implicit models presented in Section \ref{sec:partial_implicit_model}. The idea is to consider only a subset of the scenario paths of the full scenario tree and therefore obtain a solution to the implicit optimization model faster. The resulting setup plan then must be seen as a heuristic solution. The results of the partial model can be compared to the PH algorithm in a fair manner, by simply considering the same selected scenario paths in the PH approach.

In Table \ref{tab:PH_avg_partial_mode_128_3}, we report the results for sampling 128 scenario paths from the full scenario tree with $|\Omega|=3$ and therefore $|\Phi| = 2187$. We compare the PH algorithm with $\lambda=1$ and averaging consensus calculation to the partial implicit model. Both approaches include the same 128 sampled scenario paths. Results are reported for the 96 test instances. For 87 out of 96 instances, PH converged naturally, while for 9 instances a cyclic pattern was recognized and disrupted by significantly increasing $\rho$ by factor 10 to force convergence. 
After obtaining the solutions of both approaches -- PH and the partial model -- we evaluate them on the full scenario tree and report the resulting average cost, as well as the average runtimes in seconds. The column $\Delta^* (\%)$ reports the average cost delta in percent of the cost obtained by PH to the cost obtained by solution of the partial model, when both are evaluated on the full scenario tree. The cost comparison of PH in this section is therefore to the partial model, rather than the full model solution, since for the considered scenario trees the full implicit model cannot be solved to optimality within the time limit for a majority of the instances.

\begin{table}
    \resizebox{\textwidth}{!}{
    \centering
    \begin{tabular}{c|c|r|r|r|r|r}
    \hline
    \multicolumn{2}{c}{}&\multicolumn{2}{|c}{PH ($\lambda=1$, avg)}&\multicolumn{2}{|c|}{Partial model}& \\
    Util. & External demand & Runtime & Costs & Runtime & Costs & $\Delta^*$ (\%)\\
    \hline
        \multirow{2}{*}{50\%} 
         & end item demand  &  486.8 & 32361.7 & 734.5 & 32290.4 & 0.21\\
         & incl. component demand  & 778.6 & 34583.4 & 323.8 & 34138.8 & 1.32\\
         \hline
         \multirow{2}{*}{90\%} 
         & end item demand  & 967.5 & 67529.4 & 3699.3 & 66052.1 & 2.01\\
         & incl. component demand & 1593.9 & 71872.1 &  5012.3 & 71483.9 & 0.57\\
         \hline
    \end{tabular}
    }
    \caption{Averaging consensus: comparison of the partial implicit model considering 128 scenario paths sampled from the full scenario tree with $|\Omega| = 3$ and therefore $|\Phi|=2187$ to the PH algorithm with $\lambda = 1$ and averaging as consensus calculation considering the same 128 scenario paths. Results are reported for the 96 test instances and 
    grouped by utilization and demand type. Solutions of the PH and partial model are evaluated on the full scenario tree and the respective averaged costs are reported alongside averaged runtimes in seconds. The column $\Delta^*$ reports the average cost delta in percent of the PH solution to the partial model solution, when evaluated on the full scenario tree.}
    \label{tab:PH_avg_partial_mode_128_3}
\end{table}

For a resource utilization of 50\% and the case of end item demands only, we observe that PH results in a very small average cost delta of 0.21\% compared to the cost obtained by the partial model. Also, the average runtime of 486.8 seconds for PH is smaller than for the partial implicit model, which takes an average of 734.5 seconds to solve. For the case including component demands, PH has a longer runtime than the partial model, while performing 1.32\% worse in terms of costs. The case of low utilization combined with external demand for components seems to be better suited for the partial model. For the more challenging case of 90\% utilization, PH outperforms the partial model. For the case containing only end item demand, PH results in an average cost delta of 2.01\% to the cost obtained by the partial model. However, it only requires a fourth of the runtime of the model on average. For the case including external demand for components the implicit model leads to a small average cost delta of 0.57\%, while on average taking a third of the runtime of the partial model.

In Table \ref{tab:PH_mv_partial_mode_128_3}, we report the respective results for PH on the 96 test instances using the majority voting consensus calculation. 
As already indicated on the small scenario trees with $|\Omega|=2$, the obtained cost deltas to the implicit model are larger for the aggressive majority voting consensus calculation strategy, when compared to the averaging strategy. We obtained average cost deltas between 3.62\% and 7.45\%. Even though the runtimes in the majority voting setting are smaller than the runtimes in the averaging case -- especially for the difficult 90\% resource utilization case -- the significantly smaller cost deltas of averaging make it the superior choice.

\begin{table}
    \resizebox{\textwidth}{!}{
    \centering
    \begin{tabular}{c|c|r|r|r|r|r}
    \hline
    \multicolumn{2}{c}{}&\multicolumn{2}{|c}{PH ($\lambda=1$, mv)}&\multicolumn{2}{|c|}{Partial model}& \\
    Util. & External demand & Runtime & Costs & Runtime & Costs & $\Delta^*$ (\%)\\
    \hline
        \multirow{2}{*}{50\%} 
         & end item demand  & 464.1 & 33735.6 &  734.5 & 32290.4 & 3.62\\
         & incl. component demand  &  555.6 & 36439.3 & 323.8 & 34138.8 & 6.34\\
         \hline
         \multirow{2}{*}{90\%} 
         & end item demand  & 556.5 & 73362.7 & 3699.3 & 66052.1 & 7.45\\
         & incl. component demand & 995.6 & 74452.6 & 5012.3 & 71483.9 & 4.03\\
         \hline
    \end{tabular}
    }
    \caption{Majority voting consensus: comparison of the partial implicit model considering 128 scenario paths sampled from the full scenario tree with $|\Omega| = 3$ and therefore $|\Phi|=2187$ to the PH algorithm with $\lambda = 1$ and majority voting as consensus calculation considering the same 128 scenario paths. Results are reported for the 96 test instances and 
    grouped by utilization and demand type. Solutions of the PH and partial model are evaluated on the full scenario tree and the respective averaged costs are reported alongside averaged runtimes in seconds. The column $\Delta^*$ reports the average cost delta in percent of the PH solution to the partial model solution, when evaluated on the full scenario tree.}
    \label{tab:PH_mv_partial_mode_128_3}
\end{table}

We further experiment with scenario trees with $|\Omega| = 4$ and therefore $\Phi = 16,384$. We again use a partial implicit model with 128 scenario paths sampled from the full scenario tree. Note that, even though the number of paths stays constant for the partial model, the complexity of the model is still increasing. This is because the probability of two paths being indistinguishable up to some period $t$ decreases as the total number of scenario paths increases. Therefore also the number of decision variables effectively increases in the number of total scenario paths $|\Phi|$ from which we sample the 128 paths. Consider the case of $|\Omega| = 2$ with a total of $|\Phi| = 128$ scenario paths. The probability of two randomly selected paths being indistinguishable up to period 1 is 50\%, while up to period 2 it is 25\%. In the case of $|\Omega| = 4$ this probability reduces from 50\% to only 25\% for period 1 and from 25\% to 6.25\% for period 2. Therefore, in the case of $|\Omega| = 4$, out of the 128 sampled paths fewer paths will share nodes with other paths than in the case of $|\Omega| = 3$. This leads to a higher number of decision variables and therefore also a more complex model when $|\Omega|$ increases. Since the PH algorithm decomposes the full problem by scenario path anyways, the computational burden of PH stays more or less constant, as long as the number of sampled scenario paths stays constant, independently of the size of the full scenario tree. When the scenario paths share fewer nodes with each other, the consensus search might be more difficult, which can possibly increase the runtime of the PH algorithm. 

In Table \ref{tab:PH_avg_partial_mode_128_4}, we report the results for sampling 128 scenario paths from the full scenario tree with $|\Omega|=4$ and therefore $|\Phi| = 16,384$. For PH we consider the averaging consensus calculation. We again display average runtimes for the PH algorithm and average solution times to solve the partial implicit model. Further, we again report average obtained cost, when evaluating the PH solution, as well as the partial implicit model on the full scenario tree. The column $\Delta^* (\%)$ shows the average cost delta of the cost obtained by the PH solution to the cost obtained by the partial implicit model, when the respective solutions are evaluated on the full scenario tree. Results are reported for the 96 test instances. For 88 out of 96 instances, PH converged naturally, while for 8 instances a cyclic pattern was recognized and disrupted by significantly increasing $\rho$ by factor 10 to force convergence.

A resource utilization of 50\% and external demands for components is the only case, where the partial implicit model can on average be solved faster than the PH algorithm converges. It is also the only case where an evaluation of the obtained setup plans on the full scenario tree with $|\Omega| = 4$ is difficult, because resolving the model cannot be done within the time limit of three hours for all instances. For this case we report the differences in the best found solutions and a resulting average cost delta of 2.82\%, but note that due to the time limit this is only a reference value. 

For the remaining settings, however, PH significantly outperforms the implicit model in terms of runtime. Also the evaluation on the large scenario tree is possible within the time limit. In the case of 50\% utilization and external demand exclusively for end items, PH can on average be solved in 465.0 seconds, compared to 1884.6 seconds, which are on average needed to solve the respective partial implicit model. The average delta of the cost obtained by PH compared to the cost of the partial model in this case is 0.06\%. A slightly higher quality delta of 2.81\% is obtained for the case of 90\% resource utilization and end item demands only. The average runtime for the PH algorithm to generate this solution, however, is only fifteen minutes, compared to almost one hour for the implicit model. Finally, for the most difficult case including external demand for components, as well as a utilization of 90\%, the implicit model on average needs 5754.2 seconds, while the PH algorithm converges within 1339.2 seconds on average. The obtained solution on average lies within a 0.93\% delta to the costs obtained by the solution found by the partial implicit model.

\begin{table}
    \resizebox{\textwidth}{!}{
    \centering
    \begin{tabular}{c|c|r|r|r|r|r}
    \hline
    \multicolumn{2}{c}{}&\multicolumn{2}{|c}{PH ($\lambda=1$, avg)}&\multicolumn{2}{|c|}{Partial model}& \\
    Util. & External demand & Runtime & Costs & Runtime & Costs & $\Delta^*$ (\%)\\
    \hline
        \multirow{2}{*}{50\%} 
         & end item demand  &  465.0 & 34162.4 & 1884.6 & 34154.9 & 0.06 \\
         & incl. component demand  & 996.4 & $36917.9^\dagger$ & 477.8 & $35914.9^\dagger$ & $2.82^\dagger$\\
         \hline
         \multirow{2}{*}{90\%} 
         & end item demand  & 980.8 & 72797.7 & 3572.9 & 69718.8 & 2.81\\
         & incl. component demand & 1339.2 & 71590.1 & 5754.2 & 70895.0 & 0.93\\
         \hline
         \multicolumn{7}{l}{$^\dagger$ Avg best found solution after 3 hours time limit when resolving full scenario tree} \\
         \hline
    \end{tabular}
    }
    \caption{Averaging consensus: comparison of the partial implicit model considering 128 scenario paths sampled from the full scenario tree with $|\Omega| = 4$ and therefore $|\Phi|=16,384$ to the PH algorithm with $\lambda = 1$ and averaging as consensus calculation considering the same 128 scenario paths. Results are reported for the 96 test instances and 
    grouped by utilization and demand type. Solutions of the PH and partial model are evaluated on the full scenario tree and the respective averaged costs are reported alongside averaged runtimes in seconds. The column $\Delta^*$ reports the average cost delta in percent of the PH solution to the partial model solution, when evaluated on the full scenario tree.}
    \label{tab:PH_avg_partial_mode_128_4}
\end{table}

\section{Conclusion}
\label{sec:conclusion}
In this work, we have presented a PH algorithm with metaheuristic adjustment strategies for a multi-stage stochastic lot sizing problem with setup carry-overs. As a baseline, we have proposed a multi-stage stochastic programming model for the multi-item multi-echelon capacitated lot sizing problem under uncertain customer demand represented by discrete scenario trees. We have enhanced the implementation of this model by implicitly considering non-anticipativity constraints and significantly reducing the number of decision variables. 

For the PH approach, we have first implemented a base version of the algorithm, making use of a penalty parameter $\rho$ to guide the algorithm towards a consensus by penalizing consensus violations of the subproblem solutions. We have then adapted metaheuristic adjustment strategies for the PH algorithm that, on the one hand, modify subproblem cost parameters and, on the other hand, adjust the penalty parameter $\rho$ in order to guide the algorithm more efficiently to a consensus among the setup decisions. We have also proposed an alternative approach to calculate the implementable solution from the scenario subproblem solutions by means of a majority voting.

In our computational study, we have shown that the base version of the PH algorithm does not converge within three hours for a majority of the test instances, when small penalty multiplier $\lambda$-values are chosen. For larger $\lambda$-values it converges in reasonable time. However, the obtained solutions have a significant gap to the optimal costs. When evaluating the PH version including the metaheuristic adjustment strategies, we have seen a strong improvement of the algorithm convergence behaviour. Not only did the PH converge naturally for almost all of the tested $\lambda$-values -- for a small $\lambda$-value of 1, the obtained results were also within a 1\% range of the optimal costs.

For the case of the majority voting consensus calculation, we have observed a decrease of the runtime for the PH algorithm to converge, compared to the averaging consensus calculation. The obtained results, however, had slightly higher deltas to the optimal costs. The aggressive majority voting strategy highly penalizes subproblems that do no agree with the consensus. This leads to fast convergence but comes at the cost of decreased solution quality. 

Finally, we evaluate a partial implicit model only including a subset of the scenario paths, sampled from the full scenario tree. This can be seen as a heuristic solution approach, when the obtained solutions are then evaluated on the full scenario tree. The approach can be compared to PH, as long as the same subset of scenario paths is included in the algorithm. For the averaging consensus calculation we have observed that PH with the metaheuristic adjustment strategies is able to generate solutions within a 1\% range to the costs obtained by the partial implicit model in shorter runtimes. The relative reduction in runtime is more significant when the number of scenarios per tree node $|\Omega|$ increases. This is because the complexity of the implicit model increases with the size of the full scenario tree, even though the number of sampled scenario paths stays constant. The larger scenario trees lead to more isolated paths that do not share many nodes with other paths and therefore generate many decision variables in the model. The PH approach does not suffer from this effect in the same magnitude, which allows its performance to stay relatively constant in terms of runtime. We can therefore conclude that the proposed PH algorithm with metaheuristic adjustment strategies has the potential to generate very good setup plans on large scenario trees for difficult instances, containing external component demands and having high resource utilization. In such a setting one is highly limited in solving the multi-stage decision problem by means of an implicit optimization model.

\section*{Funding}
This research was funded in whole, or in part, by the Austrian Science Fund (FWF) [P 32954-G]. For the purpose of open access, the author has applied a CC BY public copyright licence to any Author Accepted Manuscript version arising from this submission.

\section*{Data availability statement}
The data that support the findings of this study are available in Zenodo at https://doi.org/10.5281/zenodo.15006105 upon reasonable request to Manuel Schlenkrich.

\bibliographystyle{apalike}
\bibliography{Literature.bib}

\end{document}